\documentclass{elsart}
\usepackage{epsfig}
\usepackage{amsmath}
 \newtheorem{Proposition}{Proposition}
 \newtheorem{Lemma}{Lemma}
 \newtheorem{Theorem}{Theorem}
 \newtheorem{Corollary}{Corollary}
 \newtheorem{definition}{Definition}
 
\usepackage{amssymb}

\begin{document}

\def\myint{\int \!\! d^{ {\scriptscriptstyle D} } \boldsymbol{x} \; }
\begin{frontmatter}

\title{Solitary waves for linearly coupled nonlinear Schr\"odinger equations with inhomogeneous coefficients}

\author{Juan Belmonte-Beitia$^{a}$, V\'{\i}ctor M. P\'erez-Garc\'{\i}a$^{a}$ and Pedro J. Torres$^b$}
\address{$^{a}$ Departamento de Matem\'aticas, E. T. S. de Ingenieros Industriales and Instituto de Matem\'atica
Aplicada a la Ciencia y la Ingenier\'{\i}a (IMACI), Universidad de
Castilla-La Mancha 13071 Ciudad Real, Spain.
 }
\address{$^b$ Departamento
de Matem\'atica Aplicada, Facultad de Ciencias,
 Universidad de Granada, Campus de Fuentenueva s/n, 18071 Granada, Spain. }

\begin{abstract}
Motivated by the study of matter
waves in Bose-Einstein condensates and coupled nonlinear optical systems, we study
a system of two coupled nonlinear Schr\"odinger equations with inhomogeneous
parameters, including a linear coupling. For that system we prove the existence of two
different kinds of homoclinic solutions to the origin describing solitary waves of physical relevance.
We use a Krasnoselskii fixed
point theorem together with a suitable compactness criterion.
\end{abstract}

\begin{keyword}
Nonlinear Schr\"{o}dinger systems, solitary waves, fixed
point theorems in cones
\end{keyword}
\end{frontmatter}

\section{Introduction}

Nonlinear Schr\"{o}dinger (NLS) equations  appear in a great array of
contexts \cite{Vazquez,Sulembook} as for example in
semiconductor electronics \cite{Soler,Soler2}, optics in nonlinear
media \cite{Kivshar}, photonics \cite{Hasegawa}, plasmas
\cite{Dodd}, fundamentation of quantum mechanics
\cite{fundamentals}, dynamics of accelerators \cite{Fedele},
mean-field theory of Bose-Einstein condensates \cite{Dalfovo} or
in  biomolecule dynamics \cite{Davidov}. In some of these fields and many others, the NLS equation
appears as an asymptotic limit for a slowly varying dispersive
wave envelope propagating in a nonlinear medium \cite{scott}.

In this paper we will be interested in systems of two coupled NLS equations in which there is a linear coupling between both components.

A first example which arises in the study of a spinor Bose-Einstein condensate composed of  two hyperfine states (for instance in many experiments, the $|F=1,m_{f}=-1\rangle$ and
$|F=2,m_{f}=1\rangle$ states of  ${}^{87}$Rb
atoms are used) and coupled by an optical or r.f. field. These systems have received a lot of attention both experimentally \cite{Matthews1,Matthews2,Inguscio1,Modugno,Inguscio2,Catani} and theoretically  (see e.g. Refs. \cite{Williams,Cirac,Sols,Deconnick,extra1,extra2,Malomed2,Nakamura,SHU} and references therein) since they represent the simplest mixture of different ultracold quantum degenerate gases and were studied immediately after the historial achievement of Bose-Einstein condensation in 1995 \cite{Dalfovo}.

To simplify the treatment yet preserving the
spatial aspects we assume the condensate to be magnetically tightly
confined along two of the transverse directions
to a single effective dimension \cite{PG98}.
In the mean-field approximation the system is described by the Gross-Pitaevskii equations
 ~\cite{Matthews1,Matthews2}
\begin{subequations}
\label{eq:6}
\begin{eqnarray}
i \frac{\partial \psi_1}{\partial t} = \left(L_1 +U_{11}
|\psi_1|^2+
U_{12} |\psi_2|^2\right)\psi_1 + \lambda \psi_2 ,
  \label{eq:6a} \\
i \frac{\partial \psi_2}{\partial t} = \left(L_2+U_{22}
|\psi_2|^2+ U_{21}|\psi_1|^2\right)\psi_2 + \lambda \psi_1,
\label{eq:6b}
\end{eqnarray}
\end{subequations}
where $L_j = -\partial^2/\partial x^2+V_j$
with $j=1,2$.  Equation (\ref{eq:6}) is  written
in dimensionless form: the spatial coordinates  and time  are measured in units of
$\ell=\sqrt{\hbar/(2m\omega)}$ and $1/\omega$, respectively, while the
energies and frequencies are measured in  units of $\hbar\omega$ and $\omega$,
$\omega$ being the trap frequency in the $(x,z)$-plane.
The dimensionless nonlinear coefficients, for the quasi-one-dimensional
condensate, are given by  $U_{ij}=U_{ji}=
2N  a_{ij}/\ell$ where $a_{ij}$ are the scattering lengths
for binary collisions, $N$ is the total number of atoms,
and  $\ell$ is the oscillator length in the $y$-direction.
The normalization of the wave function
$\psi$=col$(\psi_1,\psi_2)$ is then
$\int \psi^{\dagger}\psi d^2r=1$. In many practical situations for spinor condensates we have $U_{11} \simeq U_{12} = U_{21} \simeq U_{22}$ and it is customary to consider all of the scattering lengths to be equal.

An specific situation which arises also in applications to BEC happens when two different condensates do not coexist but are spatially separated and weakly coupled. In that situation, described for example in Refs. \cite{Li0,Li1} and references therein the model is equivalent to Eq. \ref{eq:6} but with different nonlinearities given by $U_{12} = U_{21} = 0$ (i.e. there is no nonlinear coupling between both equations).

Although the dynamics of coupled Bose-Einstein condensates has attracted a lot of attention in the last years, Equations (\ref{eq:6}) also arise in other physical contexts, most notably in Nonlinear Optical models, where they describe coupled optical fibers
 where the functions $u_{1}$ and $u_{2}$ describe the light field within each waveguide \cite{Malomed}. Also other physical situations in Nonlinear Optics are described by this type of model equations \cite{Kivshar}.

In this paper we will consider the physically relevant question of the existence of solitary wave solutions to models such as the one given by Eqs. (\ref{eq:6}) but in more general scenarios where all the parameters (i.e. the potential, the nonlinear coefficients and the coupling coefficients) are functions depending on the spatial variables.

The study of the properties of localized solutions and propagating waves in Bose-Einstein condensates with spatially inhomogeneous interactions has been a field of an enormous level of activity in Physics in the last few years \cite{INH1,INH2,INH3,INH4,INH5,INH6,INH7,INH8,INH9,INH10,INH11,INH12} motivated by novel experimental ways to control experimentally the interactions (i.e. optical manipulation of the Feschbach resonances \cite{exper}).

In this paper we will consider theoretically the problem of the existence of solitary wave solutions in two-component condensates when all of the physical parameters are spatially dependent, i.e. the potentials $V_j(x)$, the nonlinear coefficients $U_{ij}(x)$ and the coupling parameter $\lambda(x)$. The latter posibility has not been explored yet, neither experimentally nor theoretically. First because the analyses of spatially dependent
interactions have been focused on the case of a single component \cite{INH1,INH2,INH3,INH4,INH5,INH6,INH7,INH8,INH9,INH10,INH11,INH12}. Secondly because the coupling coefficient $\lambda$ is usually taken to be spatially homogeneous. However it is very simple to add an spatial dependence to this parameter in experiments since $\lambda$ is related to the coupling field and in the optical case this means just using a spatially dependent laser, i.e. a beam with a prescribed profile on the scale of the BEC.

Thus we will consider the physically relevant problem of studying solitary wave solutions to the problem
\begin{subequations}
\begin{eqnarray}
i \frac{\partial \psi_1}{\partial t} = \left[ -\frac{\partial^2}{\partial x^2}+a(x) \right]  \psi_1 - b(x) \psi_2 - c(x) F(\psi_1,\psi_2),
  \label{eq:6ar} \\
i \frac{\partial \psi_2}{\partial t} = \left[ -\frac{\partial^2}{\partial x^2}+d(x) \right]  \psi_2 - e(x) \psi_1 - f(x) G(\psi_1,\psi_2).
\label{eq:6br}
\end{eqnarray}
\end{subequations}
where $\psi_j$ are the complex wavefunctions defined on all $\mathbb{R}$ which for solitary waves, i.e. localized solutions, must also decay at infinity and we consider general nonlinearities $F$ and $G$ satisfying certain conditions to be made precise later. In order to look for stationary solutions we eliminate time as usual by defining
\begin{equation}
\psi_j(x,t) = u_j(x) \exp(i\omega_jt).
\end{equation}
Since all physically relevant nonlinearities are homogeneous in time and redefining the potentials $a(x),d(x)$ to incorporate the constant factors proportional to $\omega_ju_j$ we get the final equations to be studied in detail in this paper
\begin{subequations}\label{sistema}
\begin{eqnarray}
-u_{1}''(x)+a(x)u_{1}(x)-b(x)u_{2}(x)=c(x)F(u_{1},u_{2})u_{1},\\
-u_{2}''(x)+d(x)u_{2}(x)-e(x)u_{1}(x)=f(x)H(u_{1},u_{2})u_{2}.
\end{eqnarray}
\end{subequations}
i.e. a set of two coupled stationary linearly coupled nonlinear Schr\"odinger equations with spatially inhomogeneous coefficients. We will keep in mind the previous discussion from which it is clear that the two most interesting situations from the point of view of applications are
\begin{subequations}
\begin{equation}
\label{nl1}
F(u_{1},u_{2})=u_{1}^2,\quad H(u_{1},u_{2})=u_{2}^2
\end{equation}
 for weakly coupled Bose-Einstein condensates or optical fibers and
 \begin{equation}
\label{nl2}
F(u_{1},u_{2})=H(u_{1},u_{2})=u_{1}^2+u_{2}^2,
\end{equation}
\end{subequations}
for multicomponent spinor Bose-Einstein condensates  coupled by an optical field \cite{Williams}.
The
purpose of this paper is to study the existence of localized
positive solutions ${\bf u}=(u_{1},u_{2})$ of Eqs. (\ref{sistema}) satisfying
\begin{equation}\label{condiciones}
\lim_{|x|\rightarrow\infty}u_{i}(x)=0, \quad \lim_{|x|\rightarrow\infty}u'_{i}(x)=0,
\qquad i=1,2,
\end{equation}
where $a,b,c,d,e,f\in
L^{\infty}(\mathbb{R})$, in Eqs. (\ref{sistema}), are non-negative
almost everywhere and $F(u_{1},u_{2})$ and $G(u_{1},u_{2})$ are
continuous functions. Such solution ${\bf u}$ will have finite
energy, that mathematically means that it belongs to
$H^{1}(\mathbb{R})\times H^{1}(\mathbb{R})$.

Between the different solutions of Eq. (\ref{sistema}) those with
a more direct physical interest are the so called ground states,
which are positive real solutions of that system. Although we will be mostly interested on those solutions we will also consider the existence of other families of solutions which are relevant for applications.
The existence of positive solutions for different types of vector  Nonlinear Sch\"odinger Equations has attracted a lot of interest in recent years from the mathematical point of view. For instance in Refs. \cite{ACR,Ambrosetti,Eduardo2,Eduardo} the existence of positive solutions was proven using either critical point theory or variational approaches.

In that context in this paper we complement previous studies with the analysis of the richer system given by Eqs. \eqref{sistema}, which describes physical situations beyond those studied previously. To do so we will also develop a new theoretical approach to the problem based on the use of a fixed point theorem due to Krasnoselskii for completely continuous operators defined in cones of a Banach space together with a suitable study of the Green's function for the linear part of the problem. This method has been employed successfully for scalar problems on the real line \cite{Pedro} and in some  problems on bounded domains \cite{ORegan}. Using this technique, we will be able to prove the existence of positive solutions of Eqs. (\ref{sistema}) under the conditions (\ref{condiciones}). We think that, beyond its own interest, this topological approach appears to be complementary to the variational approach.

The rest of the paper is organized as follows. In Section \ref{preliminaries} some preliminary results are collected. Section \ref{ground} contains the main result about the existence of a homoclinic orbit to the origin (positive solution). In Section \ref{excited} we prove that, if in addition to the hypotheses of the main result of Section \ref{ground} we assume a set of hypotheses on the functions $a,b,c,d,e,f$ when there exists a second wave which is odd.  Finally, section \ref{branches} contains a result about the study of branches of solutions dependent on a parameter.

We use the notation $\mathbb{R^{+}}=(0,+\infty)$, $\mathbb{R_{+}}=[0,+\infty)$.For a given $a\in L^{\infty}(\mathbb{R})$, the essential infimum is denoted as $a_{*}$. The support of a given function $a$ is denoted by $\text{Supp(a)}$. The limit value of a given function $u$ in $+\infty$ (or $-\infty$) is written simply as $u(+\infty)$ (or $u(-\infty)$).

\section{Preliminaries}
\label{preliminaries}

The proof of the main results is based on a well-know fixed point theorem in cones for a completely continuous operator defined on a
Banach space, due to Krasnoselskii \cite{Krasnoselskii}. We recall the statement of this result below,
after introducing the definition of a cone (see, for example, \cite{FPT}).
\begin{definition}
Let $X$ be a Banach space and $P$ be a closed, nonempty subset of $X$. $P$ is a cone if
\begin{enumerate}
\item $\lambda x+\mu y\in P\quad \forall x,y\in P$ and $\quad\forall \lambda,\mu\in\mathbb{R}$,
\item $x,-x\in P$ implies $x=0$.
\end{enumerate}
\end{definition}
We also recall that a given operator is completely continuous if the image of a bounded set is relatively compact.

\begin{Theorem}\label{KN}
Let $X$ be a Banach space, and let $P\subset X$ be a cone in $X$. Assume $\Omega_{1}, \Omega_{2}$ are open subsets of $X$ with $0\in\Omega_{1},\overline{\Omega}_{1}\subset\Omega_{2}$ and let $T: P\cap(\overline{\Omega}_{2}\backslash\Omega_{1})\rightarrow P$ be a completely continuous operator such that one of the following conditions is satisfied
\begin{enumerate}
\item $\|Tu\|\leqslant\|u\|$, if $u\in P\cap\partial\Omega_{1}$, and $\|Tu\|\geqslant\|u\|$, if $u\in P\cap\partial\Omega_{2}$.
\item $\|Tu\|\geqslant\|u\|$, if $u\in P\cap\partial\Omega_{1}$, and $\|Tu\|\leqslant\|u\|$, if $u\in P\cap\partial\Omega_{2}$.
\end{enumerate}
Then, $T$ has at least one fixed point in $P\cap(\overline{\Omega}_{2}\backslash\Omega_{1})$.
\end{Theorem}

This result has been extensively employed in the study of
nonlinear equations \cite{Pedro,Zima} and also in the study of
boundary value nonlinear systems \cite{ORegan,periodic}. However,
for problems defined in non-compact intervals such as ours, there
is the difficulty that the Ascoli-Arzela theorem is not sufficient
for proving the complete continuity of the operator. We shall
employ the following compactness criterion (reminiscent of a
result by M. Zima \cite{Zima}) to show that the operator is
completely continuous.
\begin{Proposition}\label{Zima}
Let $\Omega\subset BC(\mathbb{R})$. Let us assume that the
functions $u\in\Omega$ are equicontinuous in each compact interval
of $\mathbb{R}$ and that for all $u\in\Omega$ we have
\begin{equation}\label{unif}
|u(x)|\leq\xi(x), \qquad \forall x\in\mathbb{R}
\end{equation}
where $\xi\in BC(\mathbb{R})$ verifies
\begin{equation}\label{limit}
 \lim_{|x|\rightarrow +\infty}\xi(x)=0.
\end{equation}
Then, $\Omega$ is relatively compact.
\end{Proposition}
\textbf{Proof:} Given $\{u_n\}_n$ a sequence of functions of
$\Omega$, we have to prove that there exists a partial sequence
which is uniformly convergent to a certain $u$. Note that the
elements of $\Omega$ are uniformly bounded by $\|\xi\|_\infty$ and
equicontinuous on compact intervals by hypothesis, therefore the
Ascoli-Arzela theorem provides partial sequence (call it again
$\{u_n\}_n$) which is uniformly convergent to a certain $u$ on
compact intervals. Of course, $u$ satisfies also $(\ref{unif})$.
Now, we have to prove that
$$
\forall \varepsilon>0, \,\exists n_0 \mbox{ s.t. }n\geq
n_0\implies \|u_n-u\|_{\infty}<\varepsilon.
$$
By using $(\ref{limit})$, fix $k>0$ such that
$\displaystyle\max_{x\in\mathbb{R}\setminus]-k,k[}|\xi(x)|<\varepsilon/4$.
On the other hand, by using the uniform convergence on compact
intervals, there exists $n_0$ such that
$\displaystyle\max_{x\in[-k,k]}|u_n(x)-u(x)|<\varepsilon/2$
for all $n\geq n_0$. Then,
\begin{multline*}
\|u_n-u\|_{\infty}\leq \max_{x\in[-k,k]}|u_n(x)-u(x)|+\max_{x\in\mathbb{R}\setminus]-k,k[}|u_n(x)-u(x)|<\\
<\frac{\varepsilon}{2}+2\max_{x\in\mathbb{R}\setminus]-k,k[}|\xi(x)|<\varepsilon,
\end{multline*}
and the proof is finished. \rule{5pt}{5pt}

In order to apply Theorem \ref{KN}, we need some information about
the properties of the Green's function. Let $a\in
L^{\infty}(\mathbb{R})$, $a_*>0$. For the homogeneous problem
\begin{subequations}
\begin{eqnarray}
&&-\phi''+a(x)\phi=0\\
&&\phi(-\infty)=0, \phi(\infty)=0,
\end{eqnarray}
\end{subequations}
the associated Green's function is given by
\begin{equation}\label{G1}
G_1(x,s)=\left\{\begin{array}{cc}
\phi_{1}(x)\phi_{2}(s), & -\infty<x\leq s<+\infty\\
\phi_{1}(s)\phi_{2}(x), & -\infty<s\leq x<+\infty \\
\end{array}\right .
\end{equation}
where $\phi_{1}, \phi_{2}$ are solutions such that
$\phi_{1}(-\infty)=0, \phi_{2}(+\infty)=0$. Moreover, $\phi_{1},
\phi_{2}$ can be chosen as positive increasing and positive
decreasing functions respectively. For a given $h(x)\in
L^1(\mathbb{R})$, the function
$u(x)=\int_{\mathbb{R}}G(x,s)h(s)ds$ is the unique solution of the
equation $-\phi''+a(x)\phi=h(x)$ in the Sobolev space
$H^1(\mathbb{R})$. In particular, $u$ and $u'$ vanishes at
$\pm\infty$. See for instance \cite{stuart}.

Note that $\phi_{1}, \phi_{2}$ intersect in a unique point
$x_{0}$. So, we can define a function $p_1\in BC(\mathbb{R})$ by
\begin{equation}
p_1(x)=\left\{\begin{array}{cc}
1/\phi_{2}(x), & x\leq x_{0},\\
1/\phi_{1}(x), & x>x_{0}.\\
\end{array}\right .
\end{equation}
The following result was proven in \cite{Pedro}.
\begin{Proposition}
The following properties for the Green's function defined by
(\ref{G1}) hold.
\begin{description}
\item{(P1)} $G_1(x,s)>0$ for every $(x,s)\in \mathbb{R} \times \mathbb{R}$.
\item{(P2)} $G_1(x,s)\leq G_1(s,s)$ for every $(x,s)\in \mathbb{R}\times\mathbb{R}$.
\item{(P3)} Given a non-empty compact subset $P\subset\mathbb{R}$, we define
\begin{equation}\label{defm}
m_1(P)=\min(\phi_{1}(\inf P), \phi_{2}(\sup P)).
\end{equation}
Then,
\begin{equation}
G_1(x,s)\geq m(P)p(s)G_1(s,s)\quad \text{for all}\quad (x,s)\in
P\times \mathbb{R}
\end{equation}
\end{description}
\end{Proposition}

Obviously, the linear operator $L_2[u]\equiv -u''+d(x)u$ with
$d_*>0$ defines a second Green's function $G_2(x,s)$ with
analogous properties involving properly defined $p_2(x),m_2(P)$.

\section{Existence of positive bound states}
\label{ground}

From now on, we will assume that $M=\text{Supp}(b)\cup
\text{Supp}(c)\cup \text{Supp}(e)\cup \text{Supp}(f)$ is a
non-empty compact set. In order to apply Theorem \ref{KN}, we take
the Banach space $X=BC(\mathbb{R})\times BC(\mathbb{R})$ with the
norm $\|{\bf u}\|=\text{max}_{i=1,2} |u_{i}|$, for ${\bf
u}=(u_{1},u_{2})\in X$. Let us define the cone
\begin{equation}
P=\left\{{\bf u}=(u_{1},u_{2})\in X: u_{1}(x), u_{2}(x)\geq 0\quad
\forall x,\quad \min_{x\in M}u_{i}(x)\geq m_{i}p_{0}^{i}\|u_{i}\|
\right\}
\end{equation}
where $p_{0}^{1}=\text{inf}_{M}p_{1}(x)$,
$p_{0}^{2}=\text{inf}_{M}p_{2}(x)$ and the constants $m_{1}\equiv
m_{1}(M)$ and $m_{2}\equiv m_{2}(M)$ are defined by property
$(P3)$. Note that the compactness of $M$ implies that
$p_{0}^{1}>0$ and $p_{0}^{2}>0$. Moreover, as from $(P3)$ it is
easy to see that $m_{i}p_{0}^{i}<1$, for $i=1,2$. Thus, one can
easily verify that $P$ is a cone in $X$.

Let $T: P\rightarrow X$ be a map with components $(T_{1},T_{2})$
defined by
\begin{subequations}
\begin{multline}
T_{1}({\bf u})(x)=\int_{\mathbb{R}}G_{1}(x,s)\left(b(s)u_{2}(s)+c(s)F(u_{1},u_{2})u_{1}(s)\right)ds\\
=\int_{M}G_{1}(x,s)\left(b(s)u_{2}(s)+c(s)F(u_{1},u_{2})u_{2}(s)\right)ds.
\end{multline}
\begin{multline}
T_{2}({\bf u})(x)=\int_{\mathbb{R}}G_{2}(x,s)\left(e(s)u_{1}(s)+f(s)H(u_{1},u_{2})u_{2}(s)\right)ds\\
=\int_{M}G_{2}(x,s)\left(e(s)u_{1}(s)+f(s)H(u_{1},u_{2})u_{2}(s)\right)ds.
\end{multline}
\end{subequations}
A fixed point of $T$ is a solution of system $(\ref{sistema})$
which belong to $H^1(\mathbb{R})\times H^1(\mathbb{R})$, and
therefore verifies the boundary conditions $(\ref{condiciones})$.

\begin{Lemma} Let us assume that
$$
F(u_{1},u_{2}), H(u_{1},u_{2})\geqslant 0\qquad\mbox{for every
}u_{1},u_{2}\geqslant 0.
$$
Then, $T(P)\subset P$.
\end{Lemma}
\textbf{Proof:} Take ${\bf u}=(u_{1},u_{2})\in P$. The property
$(P1)$ of the Green's function together with $b_*,c_*\geq 0$
implies $T_{1}{\bf u}(x)\geq 0$ for all $x$. Let us call $x_{m}$
the point where $min_{x\in M}T_{1}{\bf u}$ is attained. Then, for
all $x\in\mathbb{R}$,
\begin{eqnarray}
T_{1}({\bf u})(x_{m})&=&T_{1}(u_{1},u_{2})(x_{m})=\int_{M}G_{1}(x_{m},s)\left(b(s)u_{2}(s)+c(s)F(u_{1},u_{2})u_{1}(s)\right)ds\nonumber\\
&\geq & m_{1}\int_{M}p(s)G_{1}(s,s)\left(b(s)u_{2}(s)+c(s)F(u_{1},u_{2})u_{1}(s)\right)ds\nonumber\\
&\geq&m_{1}p_{0}^{1}\int_{M}G_{1}(x,s)\left(b(s)u_{2}(s)+c(s)F(u_{1},u_{2})u_{1}(s)\right)ds\nonumber\\
&=&m_{1}p_{0}^{1}T_{1}(u_{1},u_{2})(x)=m_{1}p_{0}^{1}T_{1}({\bf u})(x)
\end{eqnarray}
where we have used $(P2)$ and $(P3)$. In a similar way, we can prove that $T_{2}({\bf u})(x_{m})\geq m_{2}p_{0}^{2}T_{2}({\bf u})(x)$. This completes the proof. \rule{5pt}{5pt}

\begin{Lemma}\label{CCO}
$T: P\rightarrow P$ is continuous and completely continuous.
\end{Lemma}
\textbf{Proof:} The continuity is trivial. Let us prove that the
components of $T$ are completely continuous. Let $\Omega\subset P$
be a bounded set, with $C>0$ a uniform bound for its elements. The
functions of $T_1(\Omega)$ are equicontinuous on each compact
interval (in fact, the derivative is bounded in compacts). On the
other hand, for any ${\bf u}\in \Omega$,
$$
|T_1({\bf u})(x)|\leq
C\int_{\mathbb{R}}G_{1}(x,s)b(s)ds+C\max_{\|{\bf u}\|\leq
C}F(u_{1},u_{2})\int_{\mathbb{R}}G_{1}(x,s)c(s)ds=:\xi(x).
$$
Since the supports of $b,c$ are compact, $\xi\in
BC(\mathbb{R})\cap L^1(\mathbb{R})$, therefore $T_1(\Omega)$ is
relatively compact by Proposition \ref{Zima}. The proof for $T_2$
is analogous.
 \rule{5pt}{5pt}

The following one is the main result in this section.

\begin{Theorem}\label{principal}
Let us assume the following hypotheses,
\begin{description}
\item{(i)} $a_{*}, d_{*}>0, b_*,c_*,e_*,f_*\geqslant 0$.
\item{(ii)} $M$ is a non-empty compact set.
\item{(iii)} $F(u_{1},u_{2}), H(u_{1},u_{2})\geqslant 0$ for every $u_{1},u_{2}\geqslant 0$.
\item{(iv)} There exist $r_0>0$ and $\gamma,k>0$ such that given
$0<r<r_0$,
$$F(u_{1},u_{2}),H(u_{1},u_{2})<kr^\gamma, \mbox{ for all } \|{\bf
u}\|<r$$
\item{(v)} There exist $R_0>0$ and $\delta,K>0$ such that given
$R>R_0$,
$$
F(u_{1},u_{2}),H(u_{1},u_{2})>KR^\delta, \mbox{ if } u_1\in [m_1
p_0^1 R,R]\mbox{ or } u_2\in [m_2 p_0^2 R,R].
$$
\item{(vi)} $\int_{M}G_{1}(x,s)b(s)ds,\int_{M}G_{2}(x,s)e(s)ds<1$ for all $x\in\mathbb{R}$.
\end{description}
Then, there exists a non-trivial solution ${\bf u}\in X$ of the
system (\ref{sistema})-(\ref{condiciones}).
\end{Theorem}
\textbf{Proof:} We define the open sets $\Omega_{1}$ and
$\Omega_{2}$ as the open balls in $X$ centered in the origin and
with radius $r$ and $R$, respectively, to be fixed later. Let us
take ${\bf u}\in P\cap\partial\Omega_{1}$. Thus
\begin{multline*}
\|T{\bf u}\|=\text{max}_{x\in\mathbb{R}}
\left(\int_{M}G_{1}(x,s)\left(c(s)F(u_{1},u_{2})u_{1}(s)+b(s)u_{2}(s)\right)ds,\right. \\
\left.\int_{M}G_{2}(x,s)\left(f(s)H(u_{1},u_{2})u_{2}(s)+e(s)u_{1}(s)\right)ds\right) \end{multline*}
\begin{multline}
\leq kr^{\gamma+1}\max_{x\in\mathbb{R}}\left(\int_{M}G_{1}(x,s)c(s)ds,\int_{M}G_{2}(x,s)f(s)ds\right) \\
+r\max_{x\in\mathbb{R}}\left(\int_{M}G_{1}(x,s)b(s)ds,\int_{M}G_{2}(x,s)e(s)ds\right)\leq
r=\|{ \bf u}\|
\end{multline}
for $r$ sufficiently small, where we have used the hypothesis
$(iv)$ and $(vi)$.

On the other hand, let us take ${\bf u}\in
P\cap\partial\Omega_{2}$. Thus, at least a component of ${\bf u}$,
say $u_{1}$, satisfies $\|u_{1}\|_{\infty}=R$ (the case
$\|u_{2}\|_{\infty}=R$ is similar). Then,
$m_{1}p^{1}_{0}R\leqslant u_{1}(x)\leqslant R$ for all $x\in M$.
Therefore, by using hypothesis $(v)$ we get
\begin{multline}
\|T{\bf u}\|=\max_{i=1,2, x\in \mathbb{R}}(|T{\bf u}|)_{i}\geq\max_{x\in M}\int_{M}G_{1}(x,s)\left(c(s)F(u_{1},u_{2})u_{1}(s)+b(s)u_{2}(s)\right)ds\\
\geq KR^{\delta+1}\max_{x\in M}\int_{M}G_{1}(x,s)c(s)ds\geq
R=\|{ \bf u}\|
\end{multline}
for $R$ sufficiently large. For the case $\|u_{2}\|_{\infty}=R$,
everything works in the same way.

 Now, Theorem \ref{KN} guarantees that $T$ has a fixed point ${\bf u}\in P\cap(\overline{\Omega}_{2}\backslash\Omega_{1})$.
 Thus $r\leq{\bf \|u\|}\leq R$ so it is a non-trivial solution. \rule{5pt}{5pt}

Although the hypothesis $(vi)$ appears to be quite technical, its meaning can be clarified by realizing that
 the function $\int_{M}G_{1}(x,s)b(s)ds$ (resp.
$\int_{M}G_{2}(x,s)e(s)ds$) can be interpreted as the unique solution
of the equation $-u''+a(x)u=b(x)$ (resp. $-u''+d(x)u=e(x)$) with
boundary conditions $u(-\infty)=0=u(+\infty)$. Having this idea in
mind, we can formulate the following corollary.

\begin{Corollary}\label{cor1}
The later result still holds if (vi) is replaced by
\begin{description}
\item{(vi')} $\displaystyle\quad\max_{x\in\mathbb{R}}\frac{b(x)}{a(x)}<1,\qquad
\max_{x\in\mathbb{R}}\frac{e(x)}{d(x)}<1$
\end{description}
\end{Corollary}

\textbf{Proof:} Defining $u(x)=\int_{M}G_{1}(x,s)b(s)ds$, we need
to prove that $\|u\|_{\infty}<1$. Take $x_0$ such that
$u(x_0)=\|u\|_{\infty}$, then $u$ is convex in a neighborhood of
$x_0$, so from the equation $-u''+a(x)u=b(x)$ we easily get
$$
u(x_0)\leq\frac{b(x_0)}{a(x_0)}<1.
$$
In the same way, $\int_{M}G_{2}(x,s)e(s)ds<1$ and Theorem
\ref{principal} applies. \rule{5pt}{5pt}

\section{Odd solitary waves in the equation with symmetric coefficients}
\label{excited}

The aim of this section is to extend some ideas presented in
\cite[Section 4]{Pedro} for a scalar NLS equation to the vectorial case. In the previous
section, we have obtained positive solutions of the problem under
consideration, which is the case considered in the all related papers known to the authors. In this section, we prove the
existence of a new kind of solution when the coefficients are
even. In the following result,
$p^{1}_{0}=\inf_{M\cap\mathbb{R^{+}}}p_{1}(x)$,
$p^{2}_{0}=\inf_{M\cap\mathbb{R^{+}}}p_{2}(x)$ and $m_{1}\equiv
m_{1}(M\cap\mathbb{R^{+}})$, $m_{2}\equiv
m_{2}(M\cap\mathbb{R^{+}})$.

\begin{Theorem}
Under the conditions of theorem \ref{principal}, if moreover
$a,b,c,d,e,f$ are even functions and $0\not\in M$, then there
exists an odd non-trivial solution ${\bf u}\in
H^1(\mathbb{R})\times H^1(\mathbb{R})$ of equation (\ref{sistema})
such that $r\leq\|{\bf u}\|\leq R$.
\end{Theorem}

\textbf{Proof:} The proof mimics the steps of the proof of Theorem
\ref{principal}, working now with the Green's functions for the
problem on the half-line $\mathbb{R_{+}}$, with boundary
conditions $u_i(0)=0=u_i(+\infty)$. The operator
$T:BC(\mathbb{R_{+}})\times BC(\mathbb{R_{+}})\rightarrow
BC(\mathbb{R_{+}})\times BC(\mathbb{R_{+}})$ is defined by
\begin{subequations}
\begin{eqnarray}
T_{1}({\bf u})(x)&=&\int_{\mathbb{R_+}}G_{1}(x,s)\left(b(s)u_{2}(s)+c(s)F(u_{1},u_{2})u_{1}(s)\right)ds\\
T_{2}({\bf
u})(x)&=&\int_{\mathbb{R_+}}G_{2}(x,s)\left(e(s)u_{1}(s)+f(s)H(u_{1},u_{2})u_{2}(s)\right)ds
\end{eqnarray}
\end{subequations}
The adequate cone is
\begin{multline}
P=\left\{{\bf u}\in BC(\mathbb{R_{+}})\times BC(\mathbb{R_{+}}): u_{1}(0)=u_{2}(0)=0, \right.\\
\left. u_{1}(x), u_{2}(x)\geq 0, \forall x\in\mathbb{R_{+}},\min_{x\in M\cap\mathbb{R_{+}}}u_{i}(x)\geq m_{i}p_{0}^{i}\|u_{i}\|     \right\}.
\end{multline}
Then $T(P)\subset P$ and $T$ is a continuous and completely
continuous operator, since Proposition \ref{Zima} can be applied
to functions $u$ defined only in $\mathbb{R_{+}}$ and such that
$u(0)=0$ by simply extending as the zero constant function on the
negative axis. Finally, the sets $\Omega_{1}$ and $\Omega_{2}$ are
defined again as the open balls of radius $r$ and $R$,
respectively. Everything works in the same way so the repetitive
details are omitted. In conclusion, we obtain a positive
non-trivial solution ${\bf u}\in H^{1}(\mathbb{R_{+}})\times
H^{1}(\mathbb{R_{+}})$ such that $u_{1}(0)=u_{2}(0)=0$ and the odd
extension gives the desired solution.\rule{5pt}{5pt}

Note that the assumptions $0\not\in M$ is necessary in order to
have $m_{1}\neq 0$, $m_{2}\neq 0$, which is a key point in the
proper definition of the cone.

Finally, let us note that an analogous of Corollary \ref{cor1} also holds for this case.

\section{Branches of solutions in situations of physical interest}
\label{branches}

As it was said in the introduction, there are two situations arising in physical applications described by Eqs. (\ref{nl1}) and (\ref{nl2}).
Our results can also provide some information on the localization of the solutions that can be of interest in the study of branches of solutions in systems controlled by parameters. As a basic example, we consider the system (\ref{sistema}) with the nonlinear contribution given by Eq. (\ref{nl1}), i.e.
\begin{subequations}
\begin{eqnarray}\label{sistemaconparametro}
-u_{1}''(x)+a(x)u_{1}(x)= b(x)u_{2}(x)+\lambda c(x)u_{1}^3(x)\\
-u_{2}''(x)+d(x)u_{2}(x)=e(x)u_{1}(x)+\lambda f(x)u_{2}^3(x)
\end{eqnarray}
\end{subequations}
where $\lambda>0$.
\begin{Corollary}
Let us assume the conditions of Theorem \ref{principal}. Then, for
all $\lambda>0$ there exists a positive solution ${\bf
u_{\lambda}}=(u_{1\lambda},u_{2\lambda})\in
H^{1}(\mathbb{R})\times H^{1}(\mathbb{R})$ of
(\ref{sistemaconparametro}). Moreover
\begin{equation}
\lim_{\lambda\rightarrow 0^{+}}\|{\bf u_{\lambda}}\|=+\infty,\quad
 \lim_{\lambda\rightarrow +\infty}\|{\bf u_{\lambda}}\|=0.
\end{equation}
If moreover $a,b,c,d,e,f$ are even functions and $0\not\in M$, there exists a second branch of odd solutions ${\bf \tilde{u}_{\lambda}}\in H^{1}(\mathbb{R})$ with the same limiting properties.
\end{Corollary}

\textbf{Proof:} The application of Theorem \ref{principal}
requires the existence of $r_{\lambda}, R_{\lambda}$ such that
\begin{eqnarray}
&&r^{2}_{\lambda}\leqslant\left[1- \max_{x\in M}\left(\int_{M}G_{1}(x,s)b(s)ds, \int_{M}G_{2}(x,s)e(s)ds\right)\right]\times\nonumber\\
&&\left[\lambda\max_{x\in M}\left(\int_{M}G_{1}(x,s)c(s)ds,\int_{M}G_{2}(x,s)f(s)ds\right)\right]^{-1}\nonumber\\
&&\leqslant\left[\lambda m_0p_0^1\max_{x\in
M}\int_{M}G_{1}(x,s)c(s)ds\right]^{-1}\leqslant R^{2}_{\lambda}.
\end{eqnarray}
Since $\int_{M}G_{1}(x,s)b(s)ds,\int_{M}G_{2}(x,s)e(s)ds<1$ for
all $x\in\mathbb{R}$, such $r_{\lambda}, R_{\lambda}$ exist and
can be chosen so that
\begin{equation}
\lim_{\lambda\rightarrow 0^{+}}r_{\lambda}=+\infty,\quad
\lim_{\lambda\rightarrow +\infty}R_{\lambda}=0.
\end{equation}
Hence, we obtain a branch of solutions ${\bf u_{\lambda}}$ such that $r_{\lambda}\leqslant\|{\bf u_{\lambda}}\|\leqslant R_{\lambda}$, and now a passing to the limit finishes the proof. The arguments for the branch of odd solutions are analogous.
\rule{5pt}{5pt}

Thus, we obtain a bifurcation from infinity when $\lambda\rightarrow 0^{+}$.

We can also include the linear part under the effect of the parameter $\lambda$, i.e.
\begin{subequations}
\begin{eqnarray}\label{sistemaconparametro2}
-u_{1}''(x)+a(x)u_{1}(x)=\lambda\left( b(x)u_{2}(x)+ c(x)u_{1}^3(x)\right)\\
-u_{2}''(x)+d(x)u_{2}(x)=\lambda\left(e(x)u_{1}(x)+ f(x)u_{2}^3(x)\right)
\end{eqnarray}
\end{subequations}
where, now, $\lambda\in(0,m)$ is a positive parameter, with
\begin{equation}
 m=\left[\max_{x\in M}\left(\int_{M}G_{1}(x,s)b(s)ds,\int_{M}G_{2}(x,s)e(s)ds\right)\right]^{-1}.
 \end{equation}
Then, we obtain the following result.
\begin{Corollary}
Let us assume the conditions of Theorem \ref{principal}. Then, for
all $\lambda\in(0,m)$ there exists a positive solution ${\bf
u_{\lambda}}=(u_{1\lambda},u_{2\lambda})\in
H^{1}(\mathbb{R})\times H^{1}(\mathbb{R})$ of
(\ref{sistemaconparametro2}). Moreover
\begin{equation}
\lim_{\lambda\rightarrow 0^{+}}\|{\bf u_{\lambda}}\|=+\infty
\end{equation}
If moreover $a,b,c,d,e,f$ are even functions and $0\not\in M$, there exists a second branch of odd solutions ${\bf \tilde{u}_{\lambda}}\in H^{1}(\mathbb{R})$ with the same limiting properties.
\end{Corollary}

By using the same argument as in Corollary \ref{cor1}, one get an
explicit lower bound for $m$.

\textbf{Acknowledgements}

This work has been partially supported by grants FIS2006-04190, MTM2005-03483 (Ministerio de Educaci\'on y Ciencia, Spain) and PCI08-093 (Consejer\'ia de Educaci\'on y Ciencia de la Junta de Comunidades de Castilla-La Mancha, Spain)

\end{document}